\documentclass[12pt]{amsart}
\usepackage{amsfonts}
\usepackage{amssymb}
\usepackage{amsmath}
\input xy
\xyoption {all}

\usepackage{geometry}
\usepackage{amsthm}
\usepackage{amssymb}
\usepackage{latexsym}
\usepackage[dvips]{graphics}

\newtheorem{thm}{Theorem}[section] 
\newtheorem{cor}[thm]{Corollary}
\newtheorem{defn}[thm]{Definition}
\newtheorem{example}[thm]{Example}

\newtheorem{prop}[thm]{Proposition}
\newtheorem{remark}[thm]{Remark}

\newtheorem*{varcor}{Corollary \ref{finite}}
\newtheorem*{varthm1}{Theorem \ref{alg-geo-bound}}
\newtheorem*{varthm2}{Theorem \ref{top-bound}}
\newtheorem*{ack}{Acknowledgements}

\numberwithin{equation}{section}

\newcommand{\U}{{\mathcal U}}

\newcommand{\Z}{\mathbb{Z}}

\newcommand{\Q}{\mathbb{Q}}
\newcommand{\R}{\mathbb{R}}
\newcommand{\C}{\mathbb{C}}
\newcommand{\K}{\mathbb{K}}
\newcommand{\PP}{\mathbb{P}}

\begin{document}
\date{}
\title{Higher-order Alexander invariants of plane algebraic curves}
\author[Constance Leidy]{Constance Leidy}
\address{C. Leidy: Department of Mathematics,
          University of Pennsylvania,
          209 S 33rd St., Philadelphia, PA, 19104-6395,
                 USA.}
\email {cleidy@math.upenn.edu}

\author[Laurentiu Maxim ]{Laurentiu Maxim}
\address{L. Maxim : Department of Mathematics,
          University of Illinois at Chicago,
          851 S Morgan Street, Chicago, IL 60607,
                 USA.}
\email {lmaxim@math.uic.edu}



\begin{abstract}
We define new higher-order Alexander modules $\mathcal{A}_n(C)$ and
higher-order degrees $\delta_n(C)$ which are invariants of the
algebraic planar curve $C$. These come from analyzing the module
structure of the homology of certain solvable covers of the
complement of the curve $C$. These invariants are in the spirit of
those developed by T. Cochran  in \cite{C} and S. Harvey in \cite{H}
and \cite{Har}, which were used to study knots, 3-manifolds, and
finitely presented groups, respectively. We show that for curves in
general position at infinity, the higher-order degrees are finite.
This provides new obstructions on the type of groups that can arise
as fundamental groups of complements to affine curves in general
position at infinity.
\end{abstract}

\maketitle

\section{Introduction}

The study of plane singular curves is a subject going back to the
work of Zariski, who observed that the position of singularities has
an influence on the topology of the curve, and that this phenomena
can be detected by the fundamental group of the complement. However,
the fundamental group of a plane curve complement is in general
highly non-commutative, and thus difficult to handle. It is
therefore natural to look for invariants of the fundamental group
that capture information about the topology of the curve, such as
Alexander-type invariants associated to various covering spaces of
the curve complement. By analogy with the classical theory of knots
and links in a $3$-sphere, Libgober developed invariants of the
total linking number infinite cyclic cover in \cite{Li0, Li2, Li22}
and those of the universal abelian cover in \cite{Li3, Li4}. In this
paper, we consider certain solvable covers of the curve complement,
and their associated Alexander invariants.

Using techniques developed by T. Cochran, K. Orr, and P. Teichner in
\cite{COT}, S. Harvey (in \cite{H}) and T. Cochran (in \cite{C})
defined higher-order Alexander modules and higher-order degrees
associated to 3-manifolds and knots, respectively.  Harvey (in
\cite{Har}) has constructed similar invariants associated to
finitely presented groups, in general.  Our new invariants for
planar curves are constructed in the same way as these.

Let $C$ be a reduced curve in $\C^2$, and consider $\U$, the
complement of $C$ in $\C^2$, with $G=\pi_1(\U)$. The multivariable
Alexander invariant, studied in \cite{Li3}, \cite{Li4} (but see also
\cite{DM}), is defined by considering the universal abelian covering
space of $\U$ corresponding to the map $G \to G/[G,G] \cong \Z^s$,
where $s$ is the number of irreducible components of the curve. We
continue this construction by taking iterated universal torsion-free
abelian covers of $\U$ corresponding to the map $G \to G/G^{(n+1)}_r
\equiv \Gamma_n$, where $G^{(n)}_r$ is the $n^{th}$-term in the
rational derived series of $G$, (defined in \S2 below).  We define
the higher-order Alexander modules of the plane curve complement to
be $\mathcal{A}^{\Z}_n(C)=H_1(\U;\Z\Gamma_n)$.  The following is a
corollary to Theorem \ref{alg-geo-bound}:
\begin{varcor} If $C$ is a reduced curve in $\C^2$, that is in general position at
infinity, $\mathcal{A}^{\Z}_n(C)$ is a torsion $\Z\Gamma_n$-module.
\end{varcor}

Furthermore, we consider some skew Laurent polynomial rings
$\K_n[t^{\pm1}]$, which are obtained from $\Z\Gamma_n$ by inverting
the non-zero elements of a particular subring.  The advantage of
using $\K_n[t^{\pm1}]$ coefficients instead of $\Z\Gamma_n$
coefficients is that the former is a principal ideal domain.  We
define the higher-order degree of $C$ to be
$\delta_n(C)=\text{rk}_{\K_n} H_1(\U;\K_n[t^{\pm 1}])$.

Even though, in principal, the higher-order degrees may be computed
by means of Fox free calculus (cf. \cite{H}, $\S 6$), the
calculations are tedious as they depend on a presentation of the
fundamental group of the curve complement. However, in the case that
the curve is in general position at infinity, we find a bound on the
higher-order degrees.  In particular, we have the following theorem:

\begin{varthm1}Suppose $C$ is a degree $d$ curve in $\C^2$, such that its projective
completion $\bar{C}$ is transverse to the line at infinity. If $C$
has singularities $c_k$, $1 \leq k \leq l$, then
\begin{equation}\label{bound} \delta_n(C) \leq \Sigma_{k=1}^l \left(\mu(C,c_k) + 2 n_k\right)
+ 2g + d - l, \end{equation} where $\mu(C,c_k)$ is the Milnor number
associated to the singularity germ at $c_k$, $n_k$ is the number of
branches through the singularity $c_k$, and $g$ is the genus of the
normalized curve.
\end{varthm1}

We also can find a bound on the higher-order degrees of the curve in
terms of ``local'' degrees, $\bar{\delta_n^k}$, for each singularity
$c_k$ of $C$. The latter were defined and studied by Harvey
(\cite{H}).

\begin{varthm2} If $C$ satisfies the assumptions of the previous theorem, then
$$\delta_n(C) \leq \Sigma_{k=1}^l (\bar{\delta}_n^k + 2 n_k)+ 2g + d,$$
where $\bar{\delta}_n^k=\bar{\delta}_n(X_k)$ is Harvey's invariant
of the link complement $X_k$ associated to the singularity $c_k$.
\end{varthm2}

We view Theorem \ref{top-bound} as an analogue of the divisibility
properties for the infinite cyclic Alexander polynomial of the
complement as shown in \cite{Li2}.

For irreducible curves, regardless of the position of the line at
infinity, the higher-order degrees are finite and thus the
higher-order Alexander modules are torsion. However, if the line at
infinity is not transverse to the irreducible curve $C$, then the
upper bounds mentioned above will also include the contribution of
the singular points at infinity (similar to \cite{Li}, Thm. 4.3).

To complete the analogy with the case of Alexander polynomials of
the infinite cyclic cover of the complement, we also provide an
upper bound on $\delta_n(C)$ by the corresponding higher-order
Alexander invariant of the link at infinity (see
Theorem~\ref{infinity}). For a curve of degree $d$, in general
position at infinity, this is an uniform bound equal to $d(d-2)$.

\begin{ack}
The authors would like to thank Tim Cochran, Shelly Harvey, Anatoly
Libgober and Julius Shaneson for many helpful conversations about
this project.
\end{ack}

\section{Rational derived series of a group; PTFA groups}
In this section, we review the definitions and basic constructions
that we will need from \cite{H} and \cite{C}.  More details can be
found in these sources.

We begin by recalling the definition of the rational derived series
$\{G_r^{(i)}\}$ associated to any group $G$.

\begin{defn} Let $G_r ^{(0)}=G$. For $n \geq 1$, define the $n^{th}$ term of the rational
derived series of $G$ by:
$$G_r ^{(n)}=\{g \in G_r ^{(n-1)} | g^k \in
[G_r ^{(n-1)},G_r ^{(n-1)}], \ \text{for some} \ k \in \Z -\{0\}
\}.$$
\end{defn}

We denote by $\Gamma_n$ the quotient $G/G_r^{(n+1)}$. Since
$G_r^{(n)}$ is a normal subgroup of $G_r^{(i)}$ for $0 \leq i \leq
n$ (\cite{H}, Lemma 3.2), it follows that $\Gamma_n$ is a group.

The use of rational derived series, as opposed to the usual derived
series $\{G^{(n)}\}$, is necessary to avoid zero divisors in the
group ring $\Z \Gamma_n$. However, if $G$ is a knot group or a free
group, the rational derived series and the derived series coincide
(\cite{H}, p. 902). If $G$ is a finite group then $G_r^{(n)}=G$,
hence $\Gamma_n=\{1\}$ for all $n \geq 0$.

The rational derived series is defined in such a way that the
successive quotients $G_r^{(n)}/G_r^{(n+1)}$ are $\Z$-torsion-free
and abelian. In fact (\cite{H}, Lemma 3.5):
\begin{equation}\label{quot}
G_r^{(n)}/G_r^{(n+1)} \cong \left(G_r^{(n)}/[G_r^{(n)}, G_r^{(n)}]
\right)/\{\Z-\text{torsion}\}. \end{equation} If $G=\pi_1(X)$ this
says that $G_r^{(n)}/G_r^{(n+1)} \cong
H_1(X_{\Gamma_{n-1}})/\{\Z-\text{torsion}\}$, where
$X_{\Gamma_{n-1}}$ is the regular $\Gamma_{n-1}$-cover of $X$. In
particular, $G/G_r^{(1)}=G_r^{(0)}/G_r^{(1)} \cong
H_1(X)/\{\Z-\text{torsion}\} \cong \Z ^{\beta_1(X)}$.

\begin{remark}\rm \label{rm1} If $G$ is the fundamental group of a link
complement in $S^3$ or that of a plane curve complement, then
$G_r^{(1)}=G^{(1)}$ (since there is no torsion in the first homology
of the complement).
\end{remark}

\begin{defn}
A group $\Gamma$ is poly-torsion-free-abelian (PTFA) if it admits
a normal series of subgroups such that each of the successive
quotients of the series is torsion-free abelian.
\end{defn}

\begin{remark} \label{rm3}\rm We collect the following facts about PTFA
groups:

(1) Any subgroup of a PTFA group is a PTFA group.

(2) If $\Gamma$ is PTFA, then $\Z\Gamma$ is a right (and left) Ore
domain (i.e., has no zero divisors and $\Z\Gamma - \{0\}$ is a right
divisor set). Thus it embeds in its classical right ring of
quotients $\mathcal{K}$, a skew field (\cite{C}, Prop. 3.2).

(3) If $R$ is an Ore domain and $S$ is a right divisor set, then
$RS^{-1}$ is flat as a left $R$-module.  In particular,
$\mathcal{K}$ is a flat $\Z\Gamma$-module (\cite{Ste}, Prop.
II.3.5).

(4) Every module over $\mathcal{K}$ is a free module (\cite{Ste},
Prop. I.2.3). Such modules have a well-defined rank
$\text{rk}_{\mathcal{K}}$ which is additive on short exact
sequences.
\end{remark}

If $A$ is a module over an Ore domain $R$, then the rank of $A$ is
defined as: $rk(A)=\text{rk}_{\mathcal{K}}(A \otimes_R
\mathcal{K})$. In particular, $A$ is a torsion $R$-module if and
only if $A \otimes_R \mathcal{K}=0$, where $\mathcal{K}$ is the
quotient field of $R$.

\begin{prop}(\cite{H}, Cor. 3.6)
For any group $G$, $\Gamma_n=G/G_r^{(n+1)}$ is PTFA.
Thus it embeds in its classical right ring of quotients,
$\mathcal{K}_n$.
\end{prop}

Suppose $X$ is a topological space that has the homotopy type of a
connected CW-complex. Let $\Gamma$ be any group and $\phi:\pi_1(X)
\to \Gamma$ be a homomorphism. We denote by $X_{\Gamma}$ the regular
$\Gamma$-cover of $X$ associated to $\phi$. If $\phi$ is surjective,
this is the covering space associated to $\ker{\phi}$.  (For details
about the case where $\phi$ is not surjective, we refer the reader
to \S3 of \cite{C}.) Let $C(X_{\Gamma};\Z)$ be the $\Z\Gamma$-free
cellular chain complex for $X_{\Gamma}$ obtained by lifting the cell
structure of $X$. If $\mathcal{M}$ is a $\Z\Gamma$-bimodule, then
define:
$$H_{\ast}(X;\mathcal{M})=H_{\ast}(C(X_{\Gamma};\Z) \otimes_{\Z\Gamma} \mathcal{M})$$
as a right $\Z\Gamma$-module.

\begin{prop}(\cite{COT}, Prop. 2.9)\label{zero} Let $X$ be a connected CW-complex and $\Gamma$ a
PTFA group. If $\phi : \pi_1(X) \to \Gamma$ is a non-trivial
coefficient system, then $H_0(X;\Z\Gamma)$ is a torsion
$\Z\Gamma$-module.
\end{prop}

\begin{prop} (\cite{C}, Prop. 3.10) \label{torsion}
Let $X$ be a connected CW-complex and $\Gamma$ be a PTFA group.
Suppose $\pi_1(X)$ is finitely generated and $\phi : \pi_1(X) \to
\Gamma$ is a non-trivial coefficient system. Then $rk\left(
H_1(X;\Z\Gamma)\right) \leq \beta_1(X) -1$. In particular, if
$\beta_1(X)=1$ then $H_1(X;\Z\Gamma)$ is a torsion
$\Z\Gamma$-module.
\end{prop}

\section{Definitions of new invariants}

Let $C$ be a reduced curve in $\C^2$, defined by the equation:
$f=f_1\cdots f_s =0$, where $f_i$ are the irreducible factors of
$f$, and let $C_i=\{f_i=0\}$ denote the irreducible components of
$C$. Embed $\C^2$ in $\C\PP^2$ by adding the plane at infinity, $H$,
and let $\bar{C}$ be the projective curve in $\C\PP^2$ defined by
the homogenization $f^h$ of $f$. We let $\bar{C}_i=\{f_i^h=0\}$,
$i=1,\cdots,s$, be the corresponding irreducible components of
$\bar{C}$.
Let $\mathcal{U}$ be the complement $\mathbb{CP}^{2} - (\bar{C} \cup
H)$. (Alternatively, $\U$ may be regarded as the complement of the
curve $C$ in the affine space $\C^{2}$.) Then $H_1 (\mathcal{U})
\cong \mathbb{Z}^s$ (\cite{Di1}, (4.1.3), (4.1.4)), generated by the
meridian loops $\gamma_i$ about the non-singular part of each
irreducible component $\bar{C}_i$, for $i=1,\cdots, s$. If
$\gamma_\infty$ denotes the meridian about the line at infinity,
then in $H_1(\mathcal{U})$ there is a relation: $\gamma_\infty +
\sum {d_i \gamma_i} = 0$, where $d_i=deg(V_i)$.

\subsection{Higher-order Alexander modules}

We let $G=\pi_1(\U)$, $\Gamma_n=G/G_r^{(n+1)}$, and $\mathcal{K}_n$
be the classical right ring of quotients of $\Z\Gamma_n$.

\begin{defn}
We define the higher-order Alexander modules of the plane curve to
be:
$$\mathcal{A}^{\Z}_n(C)=H_1(\U;\Z\Gamma_n)=H_1(\U_{\Gamma_n};\Z)$$
where $\U_{\Gamma_n}$ is the covering of $\U$ corresponding to the
subgroup $G_r^{(n+1)}$.  That is,
$\mathcal{A}^{\Z}_n(C)=G_r^{(n+1)}/[G_r^{(n+1)},G_r^{(n+1)}]$ as a
right $\Z \Gamma_n$-module.
\end{defn}

\begin{defn}
The $n^{th}$ order rank of (the complement of) $C$ is:
$$r_n(C)=\text{rk}_{\mathcal{K}_n}H_1(\U;\mathcal{K}_n)$$
\end{defn}

\begin{remark}\label{rm2} \rm
(1) Note that
$\mathcal{A}^{\Z}_0(C)=G_r^{(1)}/[G_r^{(1)},G_r^{(1)}]=G'/G''$, by
Remark \ref{rm1}. This is just the universal abelian invariant of
the complement.

(2) $\mathcal{A}^{\Z}_n(C) / \{\Z-\text{torsion}\} = G^{(n+1)}_r /
G^{(n+2)}_r$.

(3) If $C$ is irreducible, then $\beta_1(\U)=1$. By Proposition
\ref{torsion}, it follows that $\mathcal{A}^{\Z}_n(C)$ is a torsion
module.

\end{remark}

In Corollary \ref{finite}, we show that under the assumption of
transversality at infinity, the module $\mathcal{A}^{\Z}_n(C)$ is a
torsion $\Z\Gamma_n$-module. Therefore, $r_n(C)=0$.

Since $\mathcal{U}$ is a $2$-dimensional affine variety, it has the
homotopy type of a $2$-dimensional CW-complex. Thus the modules
$H_k(\U;\Z\Gamma_n)$ are trivial for $k
> 2$ and $H_{2}(\U;\Z\Gamma_n)$ is a torsion-free $\Z\Gamma_n$-module. Moreover, we will show
that in our setting, the rank of $H_{2}(\U;\Z\Gamma_n)$ is equal to
the Euler characteristic of the complement, $\U$.

\begin{remark}\label{trivial}\rm
Assume that the universal abelian Alexander module of the
complement is \emph{trivial}, i.e. $\mathcal{A}^{\Z}_0(C)=0$.
(Note that this is the case if $G$ is abelian or finite.) Then all
higher-order Alexander modules $\mathcal{A}^{\Z}_n(C)=0$, for $n
\geq 1$, are also trivial. Indeed, by Remark \ref{rm1},
$G'=G_r^{(1)}$ and $\mathcal{A}^{\Z}_0(C)=G'/G''$. It follows that
$G^{(n)}=G'=G_r^{(1)}$, for all $n \geq 1$. From the definition of
the rational derived series, it is now easy to see that
$G_r^{(n)}=G'$ for all $n \geq1$. Therefore $\mathcal{A}^{\Z}_n(C)
\cong G_r^{(n+1)}/[G_r^{(n+1)},G_r^{(n+1)}] \cong G'/G'' =0$, for
all $n \geq 0$.
\end{remark}

\begin{example}\label{ex}\rm

(1) If $C$ is a non-singular curve in general position at infinity,
then $\pi_1(\U) \cong \Z$ (\cite{Li}), hence abelian. By the above
remark, it follows that $\mathcal{A}^{\Z}_n(C) = 0$, for all $n \geq
0$.\newline (2) Suppose $\U$ is the complement in $\C^2$ of a union
of two lines. Then $\pi_1(\U)$ is $\Z^2$. Hence
$\mathcal{A}^{\Z}_n(C) = 0$ for all $n \geq 0$.\newline (3) If
$\bar{C}$ is a reduced curve having only nodes as its singularities
(i.e., locally at each singular point, $\bar{C}$ looks like
$x^2-y^2=0$), then it is known that $\pi_1(\C\PP^2-\bar{C})$ is
abelian (e.g., see \cite{O}), thus has trivial commutator subgroup.
Under the assumption that the line at infinity is generic, this
implies that the commutator subgroup of $\pi_1(\U)$ is also trivial
(\cite{O}, Lemma 2), so $\pi_1(\U)$ is abelian. Now from
Remark~\ref{trivial} it follows that $\mathcal{A}^{\Z}_n(C) = 0$,
for all $n \geq 0$.

\end{example}

\subsection{Localized higher-order Alexander modules}
In this section we define some skew Laurent polynomial rings
$\K_n[t^{\pm1}]$, which are obtained from $\Z\Gamma_n$ by inverting
the non-zero elements of a particular subring described below. This
construction is used in \cite{H} and \cite{C} and is described in
algebraic generality in \cite{Har}. We refer to those sources for
the background definitions.

Recall our notations: $G =\pi_1(\U)$, $\Gamma_n=G/G_r^{(n+1)}$ and
$\mathcal{K}_n$ is the classical right ring of quotients of
$\Z\Gamma_n$. Let $\psi \in H^1(G;\Z) \cong \text{Hom}_{\Z}(G,\Z)$
be the primitive class representing the linking number homomorphism
$G \overset{\psi}{\to} \Z$, $\alpha \mapsto \text{lk}(\alpha,C)$.
Since the commutator subgroup of $G$ is in the kernel of $\psi$, it
follows that $\psi$ induces a well-defined epimorphism $\bar{\psi} :
\Gamma_n \to \Z$. Let $\bar{\Gamma}_n$ be the kernel of
$\bar{\psi}$. Since $\bar{\Gamma}_n$ is a subgroup of $\Gamma_n$, by
Remark \ref{rm3}, $\bar{\Gamma}_n$ is a PTFA group. Thus
$\Z\bar{\Gamma}_n$ is an Ore domain and $S_n=\Z\bar{\Gamma}_n -
\{0\}$ is a right divisor set of $\Z\bar{\Gamma}_n$. Let
$\K_n=(\Z\bar{\Gamma}_n)S_n ^{-1}$ be the right ring of quotients of
$\Z\bar{\Gamma}_n$, and set $R_n=(\Z\Gamma_n)S_n^{-1}$.

If we \emph{choose} a $t \in \Gamma_n$ such that $\bar{\psi}(t)=1$,
this yields a splitting $\phi : \Z \to \Gamma_n$ of $\bar{\psi}$. As
in Prop. 4.5 of \cite{H}, the embedding $\Z\bar{\Gamma}_n
\hookrightarrow \K_n$ extends to an isomorphism $R_n
\overset{\cong}{\to} \K_n[t^{\pm 1}]$. (However this isomorphism
depends on the choice of splitting!) It follows that $R_n$ is a
non-commutative principal left and right ideal domain, since this is
known to be true for any skew Laurent polynomial rings with
coefficients in a skew field (\cite{C}, Prop. 4.5). Also note that
by Remark \ref{rm3}, $R_n$ is a flat left $\Z\Gamma_n$-module.

\begin{defn}
The $n^{th}$-order localized Alexander module of the curve $C$ is
defined to be $\mathcal{A}_n(C)=H_1(\U;R_n)$, viewed as a right
$R_n$-module.  If we choose a splitting $\phi$ to identify $R_n$
with $\K_n[t^{\pm 1}]$, we define
$\mathcal{A}^{\phi}_n=H_1(\U;\K_n[t^{\pm 1}])$.
\end{defn}

\begin{defn}
The $n^{th}$-order degree of $C$ is defined to be:
$$\delta_n(C)=\text{rk}_{\K_n} \mathcal{A}_n(C)$$
\end{defn}

\begin{remark}\rm For any choice of $\phi$,
$\text{rk}_{\K_n} \mathcal{A}_n(C) = \text{rk}_{\K_n}
\mathcal{A}^{\phi}_n(C)$.  So although the module
$\mathcal{A}^{\phi}_n(C)$ depends on the splitting, the rank of the
module does not.
\end{remark}

The degrees $\delta_n(C)$ are integral invariants of the fundamental
group $G$ of the complement. Indeed, we have (\cite{Har}, \S1):
\begin{equation}\label{inv}
\delta_n(C)=rk_{\K_n} \left( {G^{(n+1)}_r}/[G^{(n+1)}_r,G^{(n+1)}_r]
\otimes_{\Z\bar{\Gamma}_n} \K_n \right)
\end{equation}

Furthermore, for any choice of splitting $\phi$, since $\K_n[t^{\pm
1}]$ is a principal ideal domain, there exist some nonzero $p_i(t)
\in \K_n[t^{\pm 1}]$, $i=1,\cdots,m$, such that:
$$\mathcal{A}^{\phi}_n(C) \cong \left( \oplus_{i=1}^m \frac{\K_n[t^{\pm 1}]}{p_i(t)\K_n[t^{\pm 1}]} \right) \oplus \K_n[t^{\pm 1}]^{r_n(C)}$$
Therefore, $\delta_n(C)$ is finite if and only if one of the
equivalent statements is true:
\begin{enumerate}
\item $r_n(C)=\text{rk}_{\mathcal{K}_n}H_1(\U;\mathcal{K}_n)=0$.
\item $\mathcal{A}_n(C)$ is a torsion $R_n$-module.
\item For any $\phi$, $\mathcal{A}^{\phi}_n(C)$ is a torsion
$\K_n[t^{\pm 1}]$-module. \item $\mathcal{A}^{\Z}_n(C)$ is a torsion
$\Z\Gamma_n$-module.
\end{enumerate}
If this is the case, then $\delta_n(C)$ is the sum of the degrees
of polynomials $p_i(t)$.

\bigskip

The invariant $\delta_n(C)$ is difficult to calculate, in general.
However, the special case of weighted homogeneous affine curves is
well understood:

\begin{prop}\label{wh}
Suppose $C$ is defined by a weighted homogeneous polynomial
$f(x,y)=0$ in $\C^2$, and assume that either $n >0$ or
$\beta_1(\U)>1$. Then we have:
\begin{equation}\label{whp}
\delta_n(C)=\mu(C,0)-1, \end{equation} where $\mu(C,0)$ is the
Milnor number associated to the singularity germ at the origin.  If
$\beta_1(\U)=1$, then $\delta_0(C)=\mu(C,0)$.

\end{prop}

\begin{proof}
The key observation here is the existence of a global Milnor
fibration (\cite{Di1}, (3.1.12)): $$F=\{f=1\} \hookrightarrow
\U=\C^2-C \overset{f}{\to} \C^{\ast},$$ and the fact that $F$ is
homotopy equivalent to the infinite cyclic cover of $\U$
corresponding to the kernel of the total linking number
homomorphism, $\psi$. The $\Gamma_n$-cover of $\U$ factors through
the infinite cyclic cover corresponding to $\psi$, which is
homotopy equivalent to $F$. It follows that there is an
isomorphism of $\K_n$-modules:
$$H_{\ast}(\U;R_n) \cong H_{\ast}(F;\K_n).$$
In particular, $$\delta_n(C)=\text{rk}_{\K_n}H_{1}(F;\K_n).$$ Since
$F$ has the homotopy type of a $1$-dimensional CW complex,
$H_{2}(F;\K_n)=0$. Moreover, if either $n>0$ or $\beta_1(\U)>1$, the
coefficient system $\pi_1(F) \to \bar{\Gamma}_n$ is non-trivial.
Hence, by Proposition \ref{zero}, $H_0(F;\K_n)=0$. It follows, in
this case, that
$$\delta_n(C)=-\chi(F)=\mu(C,0)-1.$$
On the other hand, if $\beta_1(\U)=1$, then
$\text{rk}_{\K_0}H_0(F;\K_0)=\text{rk}_{\Q}H_0(F;\Q)=1$.  Hence, if
$\beta_1(\U)=1$, then $\delta_0(C)=1-\chi(F)=\mu(C,0)$.

\end{proof}

\begin{example}\label{cusp}\rm
Since $f(x)=x^3-y^2$ is a weighted homogeneous polynomial, if $C$ is
the curve defined by $f=0$, it follows from Proposition \ref{wh},
that $\delta_0(C)=2$ and $\delta_n(C)=1$ for $n>0$.
\end{example}

\begin{remark}\label{fib}\rm
Due to the existence of Milnor fibrations for hypersurface
singularity germs, formula (\ref{whp}) holds for the case of any
algebraic link, by replacing $\U$ by the link complement and
$\delta_n(C)$ by Harvey's invariant of the algebraic link. For a
more general discussion on fibered $3$-manifolds, see \cite{H},
Prop. 8.4, 8.5.
\end{remark}

\begin{remark}\rm
As noted in \cite{H}, $\S6$ and $\S 8$, the higher-order Alexander
invariants $r_n(C)$ and $\delta_n(C)$ can be computed from a
presentation of the fundamental group of the curve complement, by
means of Fox free calculus.
\end{remark}

\section{Upper bounds on the higher order degree of a curve complement}
In this section, we find upper bounds for $\delta_n(C)$.  In Theorem
\ref{alg-geo-bound}, we find an upper bound in terms of the Milnor
number of each singularity.  In Theorem \ref{top-bound}, we find an
upper bound in terms of the Harvey's invariants, $\bar{\delta}_n$,
associated to each of the singular points of $C$. This result is
analogous to the divisibility properties for the infinite cyclic
Alexander polynomial of the complement (e.g., see \cite{Li1},
\cite{Li}, \cite{Li2}, \cite{M}). As a corollary to these theorems,
we have that, if $C$ is a curve in general position at infinity,
then $\delta_n(C)$ is finite, and therefore $\mathcal{A}_n(C)$ is a
torsion $\Z\Gamma_n$-module. We also give an upper bound for
$\delta_n(C)$ in terms of the higher-order degrees of the link at
infinity.

\begin{thm}\label{alg-geo-bound} Suppose $C$ is a degree $d$ curve in $\C^2$, such that its projective
completion $\bar{C}$ is transverse to the line at infinity, $H$. If
$C$ has singularities $c_k$, $1 \leq k \leq l$, then
\begin{equation}\label{bound} \delta_n(C) \leq \Sigma_{k=1}^l \left(\mu(C,c_k) + 2 n_k\right)
+ 2g + d - l, \end{equation} where $\mu(C,c_k)$ is the Milnor number
associated to the singularity germ at $c_k$, $n_k$ is the number of
branches through the singularity $c_k$, and $g$ is the genus of the
normalized curve.
\end{thm}
Before proving Theorem \ref{alg-geo-bound}, we state an immediate
corollary.
\begin{cor}\label{finite}
If $C$ is a plane curve in general position at infinity, then
$\delta_n(C) <\infty$, i.e., $\mathcal{A}^{\Z}_n(C)$ is a torsion
$\Z\Gamma_n$-module.
\end{cor}

\begin{remark}\rm
Note that the upper bound in (\ref{bound}) is independent of $n$.
\end{remark}

\begin{remark}\rm
If $C$ is an irreducible curve, then independently of the position
of the line at infinity, we have that $\beta_1(\U)=1$. By
Proposition \ref{torsion}, it follows that $\mathcal{A}^{\Z}_n(C)$
is a torsion module. However, if the curve $C$ is not in general
position at infinity, then the upper bound on $\delta_n(C)$ also
includes the contribution of the `singularities at infinity'
(similar to Thm 4.3 of \cite{Li}).
\end{remark}

\begin{proof}
We first reduce the problem to the study of the boundary, $X$, of a
regular neighborhood of $C$ in $\C^2$. In order to do this, let
$N(\bar{C})$ be a regular neighborhood of $\bar{C}$ inside $\C\PP^2$
and note that, due to the transversality assumption, the complement
$N(\bar{C})-(\bar{C} \cup H)$ can be identified with $N(C)-C$, where
$N(C)$ is a regular neighborhood of $C$ in $\C^2$. But $N(C)-C$
deformation retracts to $X=\partial N(C)$. Now by the Lefschetz
hyperplane section theorem (\cite{Di1}, p. 25), it follows that the
inclusion map induces a group epimorphism
$$\pi_1(X) \twoheadrightarrow \pi_1(U)$$ (the argument used here
is similar to the one used in the proof of Thm. 4.3 of \cite{Li}).
It follows that $\pi_1(X_{\Gamma_n}) \twoheadrightarrow
\pi_1(U_{\Gamma_n})$. Hence, $H_1(X;\Z\Gamma_n) \twoheadrightarrow
 H_1(\U;\Z\Gamma_n)$. Since $R_n$ is a flat
$\Z\Gamma_n$-module, there is an $R_n$-module epimorphism
$H_1(X;R_n) \twoheadrightarrow H_1(\U;R_n)$. Therefore, we have:
$$\delta_n(C) = \text{rk}_{\K_n} H_1(\U;R_n) \leq \text{rk}_{\K_n} H_1(X;R_n).$$
Hence, it suffices to show that $\text{rk}_{\K_n} H_1(X;R_n)$ is
finite. This will follow by a Mayer-Vietoris sequence argument.

Let $F$ be the (abstract) surface obtained from $C$ by removing
disks $D_1 \cup \cdots \cup D_{n_k}$ around each singular point
$c_k$ of $C$. Let $N=F \times S^1$. The boundary of $N$ is a union
of disjoint tori $T^k_1 \cup \cdots \cup T^k _{n_k}$ for
$k=1,\cdots, l$, where $l$ is the number of singular points of
$C$. For each singular point $c_k$ of $C$ we let $(S^3_k, L_k)$ be
the link pair of $c_k$, and denote by $X_k$ the link exterior,
$S^3_k-L_k$. Then $X$ is obtained from $N$ by gluing the link
exteriors $X_k$ along the tori $T^k_i$ for $i=1,\cdots,n_k$:
$$X=N \cup_{\sqcup_i{T^k_i}}(\sqcup_{k=1}^l X_k).$$
The gluing map sends each longitude of $L_k$ to the restriction of a
section in $N$, and each meridian to a fiber of $N$.

We consider the Mayer-Vietoris sequence in homology associated to
the above cover of $X$ and with coefficients in $R_n$:
\[
\begin{aligned} \cdots \to \oplus_{k,i} H_1(T^k_i;R_n) \overset{\Psi}{\to} H_1(N;R_n)
\oplus \left( \oplus_{k=1}^l H_1(X_k;R_n) \right) \to H_1(X;R_n) \\
\to \oplus_{k,i} H_0(T^k_i;R_n) \to H_0(N;R_n) \oplus \left(
\oplus_{k=1}^l H_0(X_k;R_n) \right) \to H_0(X;R_n) \to 0
\end{aligned}
\]
From Remark \ref{rm3}, we have:\begin{align}\label{rankeq}
\text{rk}_{\K_n} H_1(X;R_n) = \text{rk}_{\K_n} H_1(N;R_n) +
\Sigma_{k=1}^l \text{rk}_{\K_n} H_1(X_k;R_n) - \Sigma_{k,i}
\text{rk}_{\K_n} H_1(T^k_i;R_n) \notag \\ + \text{rk}_{\K_n}
\ker(\Psi)
 + \Sigma_{k,i} \text{rk}_{\K_n}
H_0(T^k_i;R_n) - \text{rk}_{\K_n} H_0(N;R_n) \\ - \Sigma_{k=1}^l
\text{rk}_{\K_n} H_0(X_k;R_n) +  \text{rk}_{\K_n}
H_0(X;R_n).\notag\end{align}

Recall that, for each singular point $c_k$ of $C$, the coefficient
system $R_n$ on $X_k$ is induced by the following composition of
maps:
$$\Z\pi_1(X_k) \to \Z\pi_1(X) \to \Z\pi_1(\U) \to \Z\Gamma_n \to R_n.$$
Since each $X_k$ fibers over $S^1$ with Milnor fiber $F_k$, the
$\Gamma_n$-cover of $X_k$ factors through the infinite cyclic cover
of $X_k$ which is homeomorphic to $F_k \times \R$. Therefore we have
the following isomorphisms of $\K_n$-modules:
$$H_*(X_k;R_n) \cong H_*(F_k;\K_n).$$
Since $F_k$ has the homotopy type of a wedge of circles,
$H_2(F_k;\K_n)=0$. Therefore,
$$\chi(F_k)=-\text{rk}_{\K_n} H_1(F_k;\K_n) + \text{rk}_{\K_n} H_0(F_k;\K_n)
= -\text{rk}_{\K_n} H_1(X_k;R_n) + \text{rk}_{\K_n} H_0(X_k;R_n).$$
Similar, since $N=F \times S^1$, the $\Gamma_n$-cover of $N$ factors
through the infinite cyclic cover of $N$ which is homeomorphic to $F
\times \R$. So if $F_n$ denotes the corresponding $\Gamma_n$-cover
of $F$, then $F_n$ is a non-compact surface and we have
$H_2(F;\K_n)=0$. Therefore:
$$\chi(F)=-\text{rk}_{\K_n} H_1(F;\K_n) + \text{rk}_{\K_n}
H_0(F;\K_n) = -\text{rk}_{\K_n} H_1(N;R_n) + \text{rk}_{\K_n}
H_0(N;R_n).$$ Finally, for each $k$ and $i$, we have:
$$0=\chi(S^1)=-\text{rk}_{\K_n} H_1(S^1;\K_n) + \text{rk}_{\K_n}
H_0(S^1;\K_n) = -\text{rk}_{\K_n} H_1(T^k_i;R_n) + \text{rk}_{\K_n}
H_0(T^k_i;R_n).$$

Now we can rewrite equation (\ref{rankeq}) as follows:
\begin{equation}\label{rankeq2}\text{rk}_{\K_n} H_1(X;R_n) = -\Sigma_{k=1}^l \chi(F_k) -
\chi(F)
  + \text{rk}_{\K_n} \ker(\Psi) + \text{rk}_{\K_n}
  H_0(X;R_n).\end{equation}

Since $\pi_1(X) \twoheadrightarrow \pi(U) \twoheadrightarrow
\Gamma_n$ is an epimorphism, it follows that the $\Gamma_n$-cover of
$X$ is connected, thus yielding that $\text{rk}_{\K_n}H_0(X;R_n)=1$.

Since $\Psi : \oplus_{k,i} H_1(T^k_i;R_n) \to H_1(N;R_n)$, it
follows that $\text{rk}_{\K_n} \ker(\Psi) \leq \Sigma_{k,i}
\text{rk}_{\K_n} H_1(T^k_i;R_n)$.  For each $k$ and $i$, we have
that: $$\text{rk}_{\K_n} H_1(T^k_i;R_n)=\text{rk}_{\K_n}
H_0(T^k_i;R_n)=\text{rk}_{\K_n} H_0(S^1;\K_n) \leq 1,$$ since $S^1$
is connected.  Therefore, $\text{rk}_{\K_n} \ker(\Psi)$ is less than
or equal to the number of tori, which is $\Sigma_{k=1}^l
 n_k$ where $n_k$ is the number of branches through the singularity
 $c_k$.  From equation (\ref{rankeq2}) we have the following:
\begin{equation}\label{rank3}\text{rk}_{\K_n} H_1(X;R_n) \leq \Sigma_{k=1}^l (-\chi(F_k) +
n_k) - \chi(F) + 1.\end{equation}

Furthermore, $-\chi(F_k) = \mu(C,c_k)-1$ and $-\chi(F) \leq 2g +
\sum_k {n_k} + d - 1$, where $g$ is the genus of the normalized
curve and $d$ is the degree of the curve, i.e. the number of
`punctures at infinity'. It follows that:
$$\delta_n(C) \leq \text{rk}_{\K_n} H_1(X;R_n) \leq \Sigma_{k=1}^l \left(\mu(C,c_k) + 2 n_k\right) +
2g + d - l.$$

\end{proof}

\begin{thm}\label{top-bound} Suppose $C$ is a degree $d$ curve in $\C^2$, such that its projective
completion $\bar{C}$ is transverse to the line at infinity, $H$. If
$C$ has singularities $c_k$, $1 \leq k \leq l$, then
$$\delta_n(C) \leq \Sigma_{k=1}^l (\bar{\delta}_n^k + 2 n_k)+ 2g + d,$$
where $\bar{\delta}_n^k=\bar{\delta}_n(X_k)$ is Harvey's invariant
of the link complement $X_k$ associated to the singularity $c_k$,
$n_k$ is the number of branches through the singularity $c_k$, and
$g$ is the genus of the normalized curve.
\end{thm}

\begin{proof}We have equation (\ref{rank3}) in the above proof:
$$\text{rk}_{\K_n} H_1(X;R_n) \leq \Sigma_{k=1}^l (-\chi(F_k) + n_k)
- \chi(F) + 1.$$  Furthermore, $-\chi(F) \leq 2g + \sum_k {n_k} + d
- 1$. From Prop. 8.4 of \cite{H}, we have
$$\bar{\delta}_n^k=\bar{\delta}_n(X_k)=\begin{cases} -\chi(F_k) & n \neq
0 \text{ or } \beta_1(X_k) \neq 1 \\1-\chi(F_k) & n=0 \text{ and }
\beta_1(X_k)=1.\end{cases}$$ In particular, $-\chi(F_k) \leq
\bar{\delta}_n^k$, which proves the theorem.
\end{proof}

We can also give a topological estimate for the rank of the
torsion-free $\Z\Gamma_n$-module $H_2(\U;\Z\Gamma_n)$:
\begin{cor}
The rank of the torsion-free $\Z\Gamma_n$-module
$H_2(\U;\Z\Gamma_n)$ is equal to the Euler characteristic
$\chi(\U)$ of the curve complement.
\end{cor}

\begin{proof}
Let $\mathcal{C}$ be the equivariant complex $$0 \to C_2 \to C_1 \to
C_1 \to 0$$ of free $\Z\Gamma_n$-modules, obtained by lifting the
cell structure of $\U$ to $\U_{\Gamma_n}$, the $\Gamma_n$-covering
of $\U$. Then $\chi(\mathcal{C})=\chi(\U)$. On the other hand,
$$\chi(\mathcal{C})=\sum_{i=0}^2 (-1)^i
\text{rk}_{\mathcal{K}_n}H_i(\mathcal{C} \otimes_{\Z\Gamma_n}
\mathcal{K}_n) = \sum_{i=0}^2 (-1)^i rk H_i(\mathcal{C}).$$
Therefore, by Prop. \ref{zero} and Cor. \ref{finite}, it follows
that $\chi(\mathcal{C})= rk H_2(\U;\Z\Gamma_n)$ and the claim
follows.
\end{proof}

We will end this section by relating the higher-order degrees of the
curve $C$ to the higher-order degrees of the link at infinity. Let
$S^3_{\infty}$ be a sphere in $\C^2$ of a sufficiently large radius,
or equivalently, the boundary of a small tubular neighborhood of the
hyperplane $H$ at infinity. Let $C_{\infty}=S^3_{\infty} \cap C$ be
the link of $C$ at infinity, and denote by $X_{\infty}$ its
complement $S^3_{\infty} - C_{\infty}$. Let $G_{\infty}$ denote the
fundamental group of $X_{\infty}$. We define $\delta_n^{\infty}$ to
be the $\K_n$-rank of $H_1(X_{\infty}; R_n)$ , where the coefficient
system is induced via the composition of maps: $$\Z G_{\infty} \to
\Z G \to \Z\Gamma_n \to R_n.$$

We are now ready to prove the following theorem, similar in flavor
to results on the infinite cyclic and universal abelian Alexander
invariants (see \cite{Li1}, \cite{Li}, \cite{Li2}, \cite{DM},
\cite{M}):

\begin{thm}\label{infinity}
\begin{equation}\label{inf}
\delta_n(C) \leq \delta_n^{\infty}. \end{equation}
\end{thm}

\begin{proof}
We note that there is a group epimorphism $G_{\infty}
\twoheadrightarrow G$. Indeed, $X_{\infty}$ is homotopy equivalent
to $N(H)-(\bar{C} \cup H)$, where $N(H)$ is a tubular neighborhood
of $H$ in $\C\PP^2$ whose boundary is $S^3_{\infty}$. If $L$ is a
generic line in $\C\PP^2$, which can be assumed to be contained in
$N(H)$, then by the Lefschetz theorem, it follows that the
composition
$$\pi_1(L- L \cap(\bar{C} \cup H)) \to \pi_1(N(H)-(\bar{C} \cup H)) \to \pi_1
(\C\PP^2 - (\bar{C} \cup H))$$ is surjective, thus proving our
claim (this is the same argument as the one used in \cite{Li},
Thm. 4.5).

It follows that there is a $\Z\Gamma_n$-module epimorphism
$$H_1(X_{\infty};\Z\Gamma_n) \twoheadrightarrow H_1(\U;\Z\Gamma_n).$$
Since $R_n$ is a flat $\Z\Gamma_n$-module, we also get a
$R_n$-module epimorphism:
$$H_1(X_{\infty};R_n) \twoheadrightarrow H_1(\U;R_n).$$
This proves the inequality (\ref{inf}).

\end{proof}

For a curve in general position at infinity, this yields a uniform
upper bound on the higher-order degrees of the curve, which is
independent of the local type of singularities and the number of
singular points of the curve:

\begin{cor}\label{bestb}
If $C$ is a curve of degree $d$, in general position at infinity,
then:
\begin{equation}\label{deg}
\delta_n(C) \leq d(d-2) \ , \ \ \text{for all n}.
\end{equation}
\end{cor}

\begin{proof}
The claim follows by noting that if $C$ is transversal to the line
at infinity, then $C_{\infty}$ is the Hopf link on $d_1+\cdots
+d_s=d$ components (i.e., the union of fibers of the Hopf
fibration), thus an algebraic link. By the argument used in the
proof of Proposition \ref{wh}, it follows that
$\delta_n^{\infty}=\mu_{\infty}-1$, where $\mu_{\infty}$ is the
Milnor number associated to the link at infinity. On the other hand,
$\mu_{\infty}$ is the degree of the Alexander polynomial of the link
at infinity, so it is equal to $d(d-2)+1$ (cf. \cite{Li2}). The
inequality (\ref{deg}) follows now from Theorem \ref{infinity}.

\end{proof}

\section{Examples}

In this section, we calculate the higher-order degrees for some of
the classical examples of irreducible curves, including general
cuspidal curves, Zariski's sextics with $6$ cusps, Oka's curves, and
branched loci of generic projections.

We begin with the following:
\begin{prop}\label{fact} Let $C \subset \C^2$ be an irreducible affine curve.
Let $G=\pi_1(\C^2-C)$, and denote by $\Delta_C(t)$ the Alexander
polynomial of the curve complement. If $\Delta_C(t)=1$, then
$\delta_n(C)=0$ for all $n$. Moreover, in this case,
$\mathcal{A}^{\Z}_n(C) \cong \mathcal{A}^{\Z}_0(C)$ as
$\Z[G/G']$-modules, for all $n$.

\end{prop}

\begin{proof} As $C$ is an irreducible affine curve, we have $G/G' \cong
\Z$. Hence $G' \cong G'_r$. The Alexander polynomial $\Delta_C(t)$
is the order of the infinite cyclic (and universal abelian)
Alexander module of the complement, that is $G'/G'' \otimes \Q$,
regarded as a $\Q[\Z]$-module under the action of the covering
transformations group $G/G'$ (cf. \cite{Li0, Li2, Li22}). Since $C$
is irreducible, the infinite cyclic Alexander module is a torsion
$\Q[t,t^{-1}]$-module, regardless of the position of the line at
infinity (cf. \cite{Li0}). So its order $\Delta_C(t)$ is
well-defined. Moreover, $\Delta_C(t)$ can be normalized so that
$\Delta_C(1)=1$ (\cite{Li0}).

The triviality of the Alexander polynomial means that the universal
abelian module $G'/G'' \otimes \Q$ is trivial, i.e. ${G'}/{G''}$ is
a torsion abelian group. By $(\ref{quot})$, we obtain:
$${G_r'}/{G''_r} \cong ({G_r'}/[G'_r,G'_r])/{\{\Z-\text{torsion}\}} \cong
({G'}/{G''})/{\{\Z-\text{torsion}\}} \cong 0.$$ Hence $G''_r \cong
G'_r=G'$. It follows by induction that $G_r^{(n)} = G'$, for all $n
>0$. Therefore, for any $n$,
$$\mathcal{A}^{\Z}_n(C)={G_r^{(n+1)}}/[G_r^{(n+1)},G_r^{(n+1)}] \cong
G'/G''=\mathcal{A}^{\Z}_0(C).$$ Now recall that the higher-order
degrees of $C$ may be defined by (\ref{inv}):
$$\delta_n(C)=\text{dim}_{\K_n} \left(
{G^{(n+1)}_r}/[G^{(n+1)}_r,G^{(n+1)}_r] \otimes_{\Z\bar{\Gamma}_n}
\K_n \right),$$ where $\bar \Gamma_n$ is the kernel of $\Gamma_n
\overset{\bar{\psi}}{\to} \Z=G/G'$. The map $\bar {\psi}$ is induced
by the total linking number homomorphism, which in our setting is
just the abelianization map $G \to G/G'=\Z$.  It follows that for
all $n$ we have: $\Gamma_n=G/{G^{(n+1)}_r}=G/G'$, $\bar
\Gamma_n=G'/{G^{(n+1)}_r} \cong 0$, and $\K_n\cong \Q$. Therefore,
for all $n$,
$$\delta_n(C)=\text{dim}_{\Q} \left( {G'}/{G''} \otimes_{\Z} \Q
\right)=0.$$

\end{proof}

\begin{remark}\rm Note that if the commutator subgroup $G'$ is either perfect (i.e. $G'=G''$)
or a torsion group, then the Alexander polynomial of the irreducible
curve $C$ is trivial. In particular, this is the case if $G$ is
abelian. The following examples deal with each of these cases.
\end{remark}

\begin{example}\rm Let $\bar C \subset \C\PP^2$ be an irreducible
curve of degree $d$ which has $a$ cusps (these are locally defined
by the equation $x^2=y^3$) and $b$ nodes as the only singularities.
If $d > 6a+2b$, then by a result of Nori (cf. \cite{No}, but see
also \cite{Li22}), it follows that $\pi_1(\C\PP^2- \bar C)$ is
abelian. If we choose a generic line $H$ `at infinity' and set
$C=\bar C -H$, then as in \ref{ex} it follows that $\pi_1(\C^2-C)$
is also abelian. Hence all higher-order degrees of $C$ vanish.
\end{example}

\begin{prop}\rm Let $\bar C \subset \C\PP^2$ be a degree $d$ irreducible cuspidal curve,
i.e. it admits as singularities only nodes and cusps. Choose a
generic line $H \subset \C\PP^2$, and set $C:= \bar C -H$ and
$G=\pi_1(\C^2-C)$. If $d \not\equiv 0 \ (\text{mod} \ 6)$, then all
higher-order degrees of $C$ vanish.
\end{prop}

\begin{proof}
This follows from Proposition \ref{fact} combined with Libgober's
divisibility results for the Alexander polynomial of a curve
complement (see for instance \cite{Li2}, Theorem 4.1), by noting
that the Alexander polynomial of a cusp is $t^2-t+1$, that of a node
is $t-1$, and using the fact that the Alexander polynomial of an
irreducible curve $C$ can be normalized so that $\Delta_C(1)=1$;
moreover, all zeros of $\Delta_C(t)$ are roots of unity of order
$d$.
\end{proof}

Here is a more concrete example:

\begin{example} \label{3cusp} { Zariski's three-cuspidal quartic. }\rm

Let $\bar C \subset \C\PP^2$ be a quartic curve with three cusps as
its only singularities. Choose as above a generic line, $H$, and set
$C=\bar C -H$. Then the fundamental group of the affine complement
is given by:
$$G= \pi_1(\C^2 - C)= \langle a, b \ | \ aba=bab, a^2=b^2 \rangle.$$
It is easy to see (using for example a Redemeister-Shreier process,
see \cite{MKS}) that $G' \cong \Z/3\Z$. It follows by Proposition
\ref{fact} that $\delta_n(C)=0$, for all $n$. Moreover, the integral
higher Alexander modules are given by:
$\mathcal{A}^{\Z}_n(C)=\Z/3\Z$, for all $n$.

\end{example}

\begin{remark}\rm If $\bar C \subset \C\PP^2$ is an irreducible
quartic curve, but not a three-cuspidal quartic, then the
fundamental group $\pi_1(\C\PP^2- \bar C)$ is abelian (cf.
\cite{Di1}, Proposition 4.3). If $H$ is a generic line, and $C= \bar
C -H$, then by \cite{O}, Lemma 2, it follows that $\pi_1(\C^2-C)$ is
also abelian. Thus all higher-order degrees of such a curve vanish.
Based on this observation and the previous example, it follows that
the higher-order degrees of any irreducible quartic curve are all
zero.
\end{remark}

\bigskip

In what follows, we give examples of curves having (some)
non-trivial higher-order degrees. The key observation in these
examples is the fact that the higher-order degrees of an affine
curve are invariants of the fundamental group of the complement, see
(\ref{inv}).

\begin{example}\label{sex} { Sextics with six cusps. }\rm

(a) Let $\bar C \subset \C\PP^2$ be a curve of degree $6$ with $6$
cusps on a conic. Fix a generic line, $H$, and set $C=\bar C -H$.
Then $\pi_1(\C^2 - C)=\pi_1(\C\PP^2- \bar C \cup H)$ is isomorphic
to the fundamental group of the trefoil knot, and has Alexander
polynomial $t^2-t+1$ (see \cite{Li0}, \S 7). By Remark \ref{fib},
the higher-order degrees of $C$ are the same as Cochran-Harvey
higher-order degrees for the trefoil knot, i.e. $\delta_0(C)=2$, and
$\delta_n(C)=1$ for all $n>0$.

(b) Let $\bar C \subset \C\PP^2$ be a curve of degree $6$ with $6$
cusps as its only singular points, but this time we assume that the
six cusps are not on a conic. Then $\pi_1(\C\PP^2-\bar C)$ is
abelian, (isomorphic to $\Z_2 \times \Z_3$). Assuming the line $H$
as above is generic and setting $C=\bar C -H$, this implies that
$\pi_1(\C^2- C)$ is abelian as well. Therefore, $\delta_n(C)=0$ for
all $n \geq 0$.

From (a) and (b) we see that the higher-order order degrees of a
curve, at any level $n$, are also sensitive to the position of
singular points. An interesting open problem is to find Zariski
pairs that are distinguished by some $\delta_k$, but not
distinguished by any $\delta_n$ for $n < k$.

\end{example}

\begin{example} { Oka's curves. }\rm

M. Oka \cite{Ok} has constructed the curves $\bar C_{p,q} \subset
\C\PP^2$ ($p,q$ - relatively prime), with $pq$ singular points
locally defined by
$$x^p + y^q = 0,$$ such that $\pi_1(\C\PP^2- \bar C_{p,q})=\Z_p \ast
\Z_q$. In fact, the curve $\bar C_{p,q}$ is defined by the equation:
$$(x^p+y^p)^q +(y^q+z^q)^p=0.$$
Fix a generic line $H \subset \C\PP^2$, and set $C_{p,q} = \bar
C_{p,q}-H$. Then $\pi_1(\C^2-C_{p,q})=\pi_1(\C\PP^2- \bar C_{p,q}
\cup H)$ is isomorphic to the fundamental group of the torus knot of
type $(p,q)$. The associated Alexander polynomial is (see for
instance \cite{Li0}, \S 7):
$$\Delta(t)=\frac{(t^{pq}-1)(t-1)}{(t^p-1)(t^q-1)}.$$
By Remark \ref{fib} and Proposition \ref{wh}, we obtain:
$\delta_0(C_{p,q})=\text{deg} \Delta(t)=(p-1)(q-1)$, and
$\delta_n(C_{p,q})=pq-p-q$ for all $n
> 0$.
\end{example}

\begin{example} { Branching curves of generic projections. Braid groups. }\rm

Let $V_k$ be a degree $k$ non-singular surface in $\C\PP^3$ and
$\alpha : V_k \to \C\PP^2$ be a generic projection. If $\bar C_k
\subset \C\PP^2$ denotes the branching locus of $\alpha$, then $\bar
C_k$ is  an irreducible curve of degree $k(k-1)$ with
$k(k-1)(k-2)(k-3)/2$ nodes and $k(k-1)(k-2)$ cusps. In the case
$k=3$, one obtains as branching locus the six-cuspidal sextic with
all cusps on a conic.

If $C_k$ is the affine curve obtained from $\bar C_k$ by removing
the intersection with a generic line, then B. Moishezon (\cite{Mo})
showed that $\pi_1(\C^2-C_k)$ is Artin's braid group on $k$ strands,
$B_k$. The Reidemeister-Schreier process (\cite{MKS}) leads to the
explicit computation of $B_k'/{B_k''}$ (cf. \cite{GL}). For $k \geq
5$, $B_k'/{B_k''} = 0$, hence $C_k$ has a trivial Alexander
polynomial. By Proposition \ref{fact} we obtain that
$\delta_n(C_k)=0$, for all $n \geq 0$. For $k=3$, $B_3$ is the
fundamental group of the trefoil knot, so by Example \ref{sex}(a) we
obtain: $\delta_0(C_3)=2$ and $\delta_n(C_3)=1$ for all $n>0$.

The case $k=4$ requires more work. Here we will only calculate
$\delta_0$ and $\delta_1$ of the corresponding curve $C_4$. The
Alexander polynomial of $C_4$ is $t^2-t+1$ (see for example
\cite{Li2}), thus $\delta_0(C_4)=2$. A presentation for the braid
group on four strands is:
$$B_4=\langle \sigma_1, \sigma_2, \sigma_3 | \sigma_1\sigma_3 = \sigma_3\sigma_1,
\sigma_1\sigma_2\sigma_1 = \sigma_2\sigma_1\sigma_2,
\sigma_2\sigma_3\sigma_2 = \sigma_3\sigma_2\sigma_3\>\rangle.$$ By
using Reidemeister-Schreier techniques (see for instance,
\cite{MKS}), we can obtain a presentation for $B_4'$. (This was
calculated in \cite{GL}.)
$$B_4'=\langle p,q,a,b,c | pap^{-1}=b, pbp^{-1}=b^2c, qaq^{-1}=c, qbq^{-1}=c^3a^{-1}c, c=a^{-1}b \rangle,$$
where $p=\sigma_2\sigma_1^{-1}$, $q=\sigma_1\sigma_2\sigma_1^{-2}$,
$a=\sigma_3\sigma_1^{-1}$,
$b=\sigma_2\sigma_1^{-1}\sigma_3\sigma_2^{-1}$, and
$c=\sigma_1\sigma_2\sigma_1^{-2}\sigma_3\sigma_1\sigma_2^{-1}\sigma_1^{-1}$.
Then, $B_4'/B_4'' \cong \Z \oplus \Z$, generated by $p$ and $q$.
Notice that since $B_4'/B_4''$ is torsion-free, $(B_4)''_r = B_4''$.
Hence by (\ref{inv}), we have:
$$\delta_1 = \text{rk}_{\K_1}\left(B_4''/B_4'''\otimes_{\Z\bar{\Gamma}_1} \K_1\right),$$
where $\bar{\Gamma}_1=\ker(\bar{\psi}:B_4/B_4'' \to
B_4/B_4')=B_4'/B_4''$. Therefore, we must understand $B_4''/B_4'''$
as a $\Z[p^{\pm1},q^{\pm1}]$-module, and then determine the rank as
a $\Q(p,q)$-vector space.

Again using Reidemeister-Schreier techniques, we calculate a group
presentation for $B_4''$:
\begin{eqnarray*}B_4''=\langle \rho_{i,j}, \alpha_{i,j}
| \rho_{i,0}=1, \rho_{i,j}\alpha_{i+1,j}\rho_{i,j}^{-1} =
\alpha_{i,j}\alpha_{i,j+1},
\rho_{i,j}\alpha_{i+1,j+1}\rho_{i,j}^{-1} = \alpha_{i,j}
\left(\alpha_{i,j+1}\right)^2, \\
\alpha_{i,j+2} =
\left(\alpha_{i,j+1}\right)^2\alpha_{i,j}^{-1}\alpha_{i,j+1}
\rangle,\end{eqnarray*} where $\rho_{i,j}=p^i q^j p q^{-j}
p^{-(i+1)}$ and $\alpha_{i,j}=p^i q^j a q^{-j} p^{-i}$. Notice that
$p$ and $q$ act on $B_4''$ by conjugation. Furthermore,
$p*(q*\gamma)=\rho_{0,1}^{-1}(q*(p*\gamma))\rho_{0,1}$, for all
$\gamma \in B_4''$. Hence although the actions of $p$ and $q$ do not
commute in $B_4''$, they do commute in $B_4''/B_4'''$. In
particular, $B_4''/B_4'''$ is indeed a
$\Z[p^{\pm1},q^{\pm1}]$-module.

We have the following presentation for $B_4''/B_4'''$ as an abelian
group: $$B_4''/B_4'''=\langle \rho_{i,j}, \alpha_{i,j} |
\rho_{i,0}=0, \alpha_{i+1,j}=\alpha_{i,j}+\alpha_{i,j+1},
\alpha_{i+1,j+1}=\alpha_{i,j}+2\alpha_{i,j+1},\alpha_{i,j+2}=3\alpha_{i,j+1}-\alpha_{i,j}
\rangle.$$ To get a presentation as a
$\Z[p^{\pm1},q^{\pm1}]$-module, we note that: \begin{eqnarray*}
\rho_{i,j}&=&\prod_{k=1}^j (p^i q^{j-k} * \rho_{0,1}), \text{ for }
j
\geq 1, \\
\rho_{i,j}&=&\prod_{k=1}^{-j} (p^i q^{j-1+k} * \rho_{0,1}^{-1}),
\text{ for } j \leq -1, \\
\rho_{i,0}&=&0, \\
\alpha_{i,j}&=&p^i q^j * \alpha_{0,0}, \text{ for all }
i,j\in\Z.\end{eqnarray*} Therefore, as a
$\Z[p^{\pm1},q^{\pm1}]$-module, $B_4''/B_4'''$ is generated by
$\rho_{0,1}$ and $\alpha_{0,0}$. Furthermore, $\rho_{0,1}$ generates
a free submodule, while $\alpha_{0,0}$ generates a torsion
submodule. Hence the rank as a $\Q(p,q)$-vector space is 1.
Therefore, $\delta_1(C_4)=1$.

\emph{Note.} The background material on the constructions mentioned
in this example are beautifully explained in Libgober's papers
\cite{Li2} and \cite{Li22}. In particular, the latter contains a
summary of Moishezon's results \cite{Mo}.
\end{example}

\section{Concluding Remarks}
Although in geometric problems the fundamental group of complements
to projective curves plays a central role, by switching to the
affine setting (i.e. by removing also a generic line) no essential
information is lost. Indeed, the two groups are related by the
central extension $$0 \to \Z \to \pi_1(\C\PP^2-(\bar{C} \cup H)) \to
\pi_1(\C\PP^2 - \bar {C}) \to 0.$$

Our finiteness result on the higher-order degrees provides new
obstructions on the type of groups that can arise as fundamental
groups of complements to affine curves in general position at
infinity. Note that for a general group, one does not expect the
higher-order degrees $\delta_n$ to be finite. For instance, for a
free group with at least 2 generators the free ranks $r_n$ are
positive (cf. \cite{H}, Example 8.2) therefore $\delta_n$ is
infinite.

Similar obstructions were previously obtained by Libgober. For
example, from the study of the total linking number infinite cyclic
cover of the complement (\cite{Li0}, \cite{Li2}), it follows that
the Alexander polynomial of the (affine) curve is cyclotomic. More
precisely, this polynomial divides the product of the local
Alexander polynomials at the singular points, and its zeros are
roots of unity of order $d=\text{deg} (\bar C)$. Note also that the
fundamental group of the affine complement maps onto $\Z/d$ (this
can be seen from the above central extension). More obstructions
were derived by Libgober from the study of the universal abelian
cover of the affine complement (\cite{Li4}, \cite{Li3}): for
example, the support of the universal abelian module is contained in
the zero-set of the polynomial $ (\prod_{i=1}^s t_i^{d_i})-1=0$,
where $d_i$ are the degrees of the irreducible components of the
curve.

Our obstructions come from analyzing the solvable coverings
associated to the rational derived series of the fundamental group
of the affine complement. It would be interesting to understand how
the higher-order degrees are related to (or influenced by) the
invariants of the infinite cyclic or universal abelian covers of the
complement. Proposition \ref{fact} already provides such a relation.
In connection with the universal abelian cover, A. Libgober told us
that he proved the following result: if the codimension (in the
character torus) of support of the universal abelian Alexander
module is greater than $1$, then $\delta_0(C)=0$. Of course, this
assumption can only be satisfied if the curve is reducible, with at
least $2$ components, but it is an interesting problem to understand
for what type of curves this condition holds.

\providecommand{\bysame}{\leavevmode\hbox
to3em{\hrulefill}\thinspace}

\end{document}